\documentclass[11pt]{article}
\usepackage{amsmath,amssymb,theorem}
\usepackage{graphicx}
\usepackage{subfigure}
\usepackage{psfrag}
\usepackage{color}
\usepackage{enumerate}
\usepackage{epstopdf}
\newtheorem{thm}{Theorem}[section]
\newtheorem{conj}[thm]{Conjecture}
\newtheorem{lem}[thm]{Lemma}%[section]
\theorembodyfont{\rmfamily}
\def\pf{\bigskip\noindent {\bf Proof.}~~}

\def\pf{\bigskip\noindent {\bf{Proof.}}~~}

\title{Rainbow matchings in properly-colored hypergraphs}
\date{}
\author{
	Hao Huang \thanks{
		Department of Math and CS, Emory University, Atlanta, GA 30322, USA.
		Email: hao.huang@emory.edu. Research supported in part by the Collaboration Grants from the Simons Foundation.
	}
	\and
	Tong Li\thanks{
		Department of Mathematics,
		Shandong University, Jinan, China.
	}
	\and
		Guanghui Wang
		\thanks{
		Department of Mathematics,
		Shandong University, Jinan, China. Research supported by NNSF (No. 11471193, 11631014).
	}
}

\begin{document}
\maketitle
\begin{abstract}
A hypergraph $H$ is properly colored if for every vertex $v\in V(H)$, all the edges incident to $v$ have distinct colors. In this paper, we show that if $H_{1}$, \ldots, $H_{s}$ are properly-colored $k$-uniform hypergraphs on $n$ vertices, where $n\geq3k^{2}s$, and  $e(H_{i})>{{n}\choose {k}}-{{n-s+1}\choose {k}}$, then there exists a rainbow matching of size $s$, containing one edge from each $H_i$. This generalizes some previous results on the Erd\H{o}s Matching Conjecture.

\end{abstract}
\bigskip
\noindent {\textbf{Keywords}: rainbow matching, properly-colored hypergraphs

\baselineskip 18pt
\section{Introduction}
A \emph{$k$-uniform hypergraph} is a pair $H=(V, E)$, where $V=V(H)$ is a finite set of vertices, and $E=E(H)\subseteq {{V}\choose {k}}$ is a family of $k$-element subsets of $V$ called edges. A \emph{matching} in a hypergraph $H$ is a collection of vertex-disjoint edges. The \emph{size} of a matching is the number of edges in the matching. The matching number $\nu(H)$ is the maximum size of a matching in $H$. In 1965, Erd\H{o}s \cite{PE} asked to determine the maximum number of edges that could appear in a $k$-uniform $n$-vertex hypergraph $H$ with matching number $\nu(H)<s$, for given integer $s \leq\frac{n}{k}$. He conjectured that the problem has two extremal constructions. The first one is a hyper-clique consisting of all the $k$-subsets on $ks-1$ vertices. The other one is a $k$-uniform hypergraph on $n$ vertices containing all the edges intersecting a fixed set of $s-1$ vertices. Erd\H os  posed the following conjecture:

\begin{conj}[\cite{PE}]\label{conj1}
Every $k$-uniform hypergraph $H$ on $n$ vertices with matching number $\nu(H)<s\leq\frac{n}{k}$ satisfies $e(H)\leq \max\{{{ks-1}\choose {k}}, {{{n}\choose {k}}-{{n-s+1}\choose {k}}}\}$.
\end{conj}

The case $s=1$ is the classic Erd\H{o}s--Ko--Rado Theorem \cite{EKR}. The graph case ($k=2$) was verified in \cite{EG} by Erd\H{o}s and Gallai. The problem seems to be significantly harder for hypergraphs. When $k=3$, Frankl, R\"{o}dl and Ruci\'{n}ski \cite{FR} proved the conjecture for $s\leq \frac{n}{4}$. {\L}uczak and Mieczkowska \cite{LM} proved it for sufficiently large $s$. The $k=3$ case was finally settled by Frankl \cite{F2}. For general $k$, a short calculation shows that when $s \leq \frac{n}{k+1}$, we always have ${{{n}\choose {k}}-{{n-s+1}\choose {k}}}>{{ks-1}\choose {k}}$. For this range, the second construction is believed to be optimal. Erd\H{o}s \cite{PE} proved the conjecture for $n\geq n_{0}(k, s)$. Bollob\'{a}s, Daykin and Erd\H{o}s \cite{BDE} proved the conjecture for $n>2k^{3}(s-1)$. Huang, Loh and Sudakov \cite{HL} improved it to $n\geq 3k^{2}s$, which was further improved to $n\geq {3k^{2}s}/{\log k}$ by Frankl, {\L}uczak and Mieczkowska \cite{FL}. On the other hand, in an unpublished note, F\"{u}redi and Frankl proved the conjecture for $n\geq cks^2$, Frankl \cite{F} improved all the range above to $n\geq (2s-1)k-s+1$. Currently the best range is $n \ge \frac{5}{3}sk-\frac{2}{3}s$ by Frankl and Kupavskii \cite{FK}.

In this paper, we consider a generalization of Erd\H{o}s Matching Conjecture to properly-colored hypergraphs. A hypergraph $H$ is \emph{properly colored} if for every vertex $v\in V(H)$, all edges incident to $v$ are colored differently. A \emph{rainbow matching} in a properly-colored hypergraph $H$ is a collection of vertex disjoint edges with pairwise different colors. The \emph{size} of a rainbow matching is the number of edges in the matching. The \emph{rainbow matching number}, denoted by $\nu_{r}(H)$, is the maximum size of a rainbow matching in $H$.  Motivated by the Erd\H{o}s Matching Conjecture, we consider the following problem: how many edges can appear in a properly-colored $k$-uniform hypergraph $H$ such that its rainbow matching number satisfies $\nu_{r}(H)<s\leq\frac{n}{k}$? In fact, it is called Rainbow Tur\'{a}n problem and is well studied in  \cite{keevash}. Note that here if we let $H$ be rainbow, that is, every edge of $H$ receives distinct colors, then we obtain the original Erd\H os Matching Conjecture.

More generally, let $H_{1}$, \ldots, $H_{s}$ be properly-colored $k$-uniform hypergraphs on $n$ vertices, a rainbow matching of size $s$ in $H_{1}$, \ldots, $H_{s}$ is a collection of vertex disjoint edges ${e_1, \ldots, e_s}$ with pairwise different colors, where $e_{1}\in E(H_1), \ldots, e_{s}\in E(H_s)$. For simplicity, we call it an $s$-rainbow matching. Then what is the minimum $M$, such that by assuming $e(H_i) >M$ for every $i$, it guarantees the existance of an $s$-rainbow matching?

In this paper, we prove the following result, which generalizes Theorem 1.2 and Theorem 3.3 of \cite{HL}.

\begin{thm}\label{thm1}
Let $H_{1}$, \ldots, $H_{s}$ be properly-colored $k$-uniform hypergraphs on $n$ vertices. If $n\geq 3k^{2}s$ and every $e(H_{i})>{{n}\choose {k}}-{{n-s+1}\choose {k}}$, then there exists an $s$-rainbow matching in $H_{1}$, \ldots, $H_{s}$.
\end{thm}

\section{Preliminary results}

In this section, we list some preliminary results about ``rainbow'' hypergraphs, which is a special case of properly-colored hypergraphs. In the next section, we will prove our main theorem with the help of these results. A hypergraph $H$ is \emph{rainbow} if the colors of any two edges in  $E(H)$ are different. From now on, when we say an edge $e$ is \emph{disjoint} from a collection of edges, it means that not only $e$ is vertex-disjoint from those edges, but it also has a color different from the colors of all these edges. We start by the following lemma for graphs. Note that here although each $G_i$ is rainbow, a color may appear in more than one $G_i$'s.

\begin{lem}\label{lem1}
Let $G_{1}$, \ldots, $G_{s}$ be rainbow graphs on $n$ vertices. If $n\geq 5s$ and $e(G_{i})>{{n}\choose {2}}-{{n-s+1}\choose {2}}$, then there exists an $s$-rainbow matching in $G_{1}$, \ldots, $G_{s}$.
\end{lem}

\pf We do induction on $s$. The base case $s=1$ is trivial.
For every vertex $v\in V(G_{i})$ and $j\neq i$, let $G_{v}^{j}$ be the subgraph of $G_j$ induced by the vertex set $V(G_{j})\setminus\{v\}$. Since there are at most $n-1$ edges containing $v$ in $E(G_{j})$, we have $e(G_{v}^{j})\geq e(G_{j})-(n-1)>{{n}\choose {2}}-{{n-s+1}\choose {2}}-(n-1)={{n-1}\choose {2}}-{{(n-1)-(s-1)+1}\choose {2}}$. By induction, there exists an $(s-1)$-rainbow matching $\{e_j\}_{j\neq i}$ in $\{G_{v}^{j}\}_{j\neq i}$, which spans $2(s-1)$ vertices. So if some $G_i$ has a vertex $v$ with degree greater than $3(s-1)$, then there exists an edge $e$ in $G_i$ which contains $v$ and disjoint from the edges of the $(s-1)$-rainbow matching, which produces an $s$-rainbow matching. Hence we may assume that the maximum degree of each $G_i$ is at most $3(s-1)$.

 Now pick an arbitrary edge $uv$ in $G_1$. Assume the color of $uv$ is $c(uv)$. Then we delete the vertices $u$, $v$ and the edge colored by $c(uv)$ in $G_{2}, \ldots, G_{s}$. Denote the resulting graphs by $G'_{2}, \ldots, G'_{s}$. We can see that when $n\geq5s$, for each $i\in \{2, \ldots, s\}$, we have $e(G'_{i})>{{n}\choose {2}}-{{n-s+1}\choose {2}}-2\cdot3(s-1)-1>{{n-2}\choose {2}}-{{(n-2)-(s-1)+1}\choose {2}}$. By induction on $s$, there exists an $(s-1)$-rainbow matching in the graphs $G'_{2}, \ldots, G'_{s}$. Taking these $s-1$ edges with the edge $uv$, we obtain an $s$-rainbow matching in $G_{1}$, \ldots, $G_{s}$.

\hfill\vrule height3pt width6pt depth2pt

\begin{lem}\label{lem2}
Let $H_{1}$, \ldots, $H_{s}$ be rainbow $k$-uniform hypergraphs on $n$ vertices. If $n\geq 3k^{2}s$ and $e(H_{i})>{{n}\choose {k}}-{{n-s+1}\choose {k}}$, then there exists an $s$-rainbow matching in $H_{1}$, \ldots, $H_{s}$.
\end{lem}

\pf We do induction on both $k$ and $s$. According to Lemma \ref{lem1}, the case $k=2$ holds for every $s$ and $n\geq 5s$. And for every $k$, the case $s=1$ is trivial.
We first consider the situation when some $H_i$ has a vertex $v$ with degree greater than $k(s-1){{n-2}\choose {k-2}}+s-1$. For every vertex $v\in V(H_{i})$ and $j\neq i$, let $H_{v}^{j}$ be the subgraph of $H_j$ induced by the vertex set $V(H_{j})\setminus \{v\}$. Since there are at most ${n-1}\choose {k-1}$ edges containing $v$ in $E(H_{j})$, we have $e(H_{v}^{j})\geq e(H_{j})-{{n-1}\choose{k-1}}>{{n}\choose {k}}-{{n-s+1}\choose {k}}-{{n-1}\choose{k-1}}={{n-1}\choose {k}}-{{(n-1)-(s-1)+1}\choose {k}}$. By inductive hypothesis for the case $(n-1, k, s-1)$, there exists an $(s-1)$-rainbow matching $\{e_j\}_{j\neq i}$ in $\{H_{v}^{j}\}_{j\neq i}$, which spans $k(s-1)$ vertices. So if some $H_i$ has a vertex $v$ with degree greater than $k(s-1){{n-2}\choose {k-2}}+s-1$, then there exists an edge $e$ in $E(H_i)$ which contains $v$ and disjoint from  the edges of the $(s-1)$-rainbow matching, which produces an $s$-rainbow matching. Hence we may assume that the maximum degree in each hypergraph $H_i$ is at most $k(s-1){{n-2}\choose {k-2}}+s-1$.

By induction on $s$, we know that for every $i$ there exists an $(s-1)$-rainbow matching in the hypergraphs $\{H_j\}_{j\neq i}$, spanning $k(s-1)$ vertices. If for some $i$, the $s$-th largest degree of $H_i$ is at most $2(s-1){{n-2}\choose {k-2}}+s-1$, then the sum of degrees of these $k(s-1)$ vertices in $H_i$ is at most \medskip

 $(s-1)[k(s-1){{n-2}\choose {k-2}}+s-1]+(s-1)(k-1)[2(s-1){{n-2}\choose {k-2}}+s-1]=(3k-2)(s-1)^{2}{{n-2}\choose {k-2}}+(s-1)^{2}k$.\medskip

Since $n\geq 3k^{2}s$, we have $e(H_i)>{{n}\choose {k}}-{{n-s+1}\choose {k}}>(s-1)^{2}(3k-\frac{1}{2}){{n-2}\choose {k-2}}>(3k-2)(s-1)^{2}{{n-2}\choose {k-2}}+(s-1)^{2}k+s-1$, which guarantees the existence of an edge in $H_i$ which is disjoint from the previous $(s-1)$-rainbow matching in $\{H_j\}_{j\neq i}$, which produces an $s$-rainbow matching. So we may assume that each $H_i$ contains at least $s$ vertices with degree above  $2(s-1){{n-2}\choose {k-2}}+s-1$.

Now we may greedily select distinct vertices $v_{i}\in V(H_{i})$, such that for each $1\leq i\leq s$, the degree of $v_i$ in $H_i$ exceeds $2(s-1){{n-2}\choose {k-2}}+s-1$. Consider all the subsets of $V(H_{i})\setminus \{v_{1}, \ldots, v_{s}\}$ which together with $v_i$ form an edge of $H_i$. Denote the $(k-1)$-uniform hypergraph by $H_{i}'$. Then $e(H_{i}')>2(s-1){{n-2}\choose {k-2}}+s-1-(s-1){{n-2}\choose {k-2}}>{{n-s}\choose {k-1}}-{{n-2s+1}\choose {k-1}}$. By the inductive hypothesis for the case $(n-s, k-1, s)$, there exists an $s$-rainbow matching $\{e_{i}\}_{1\leq i\leq s}$ in $\{H_{i}'\}_{1\leq i\leq s}$. Taking the edges $e_{i}\bigcup \{v_i\}$, we obtain an $s$-rainbow matching in $\{H_{i}\}_{1\leq i\leq s}$.

\hfill\vrule height3pt width6pt depth2pt

%%%%%%%
%%%%%%%
%%%%%%%
%%%%%%%
%%%%%%%
\section{Main Theorem}

In this section we  prove our main result, Theorem \ref{thm1}, using induction and Lemma \ref{lem2}.

%\begin{thm}
%Let $H_{1}$,\ldots, $H_{s}$ be properly-colored $k$-uniform hypergraphs on $n$ vertices. If $n\geq 3k^{2}s$ and $e(H_{i})>{{n}\choose {k}}-{{n-s+1}\choose {k}}$, then there exists an $s$-rainbow matching in $H_{1}$, \ldots, $H_{s}$.
%\end{thm}

\pf We split our proof into two cases.

\textbf{Case 1: $k=2$.} Now $H_{1}$, \ldots, $H_{s}$ are properly-colored graphs. We do induction on $s$. The base case $s=1$ is trivial.
    For every vertex $v\in V(H_{i})$ and $j\neq i$, let $H_{v}^{j}$ be the subgraph of $H_j$ induced by the vertex set $V(H_{j})\setminus \{v\}$. Since there are at most $n-1$ edges containing $v$ in $E(H_{j})$, we have $e(H_{v}^{j})\geq e(H_{j})-(n-1)>{{n}\choose {2}}-{{n-s+1}\choose {2}}-(n-1)={{n-1}\choose {2}}-{{(n-1)-(s-1)+1}\choose {2}}$. By induction, there exists an $(s-1)$-rainbow matching $\{e_j\}_{j\neq i}$ in $\{H_{v}^{j}\}_{j\neq i}$, which spans $2(s-1)$ vertices. So if some $H_i$ has a vertex $v$ of degree greater than $3(s-1)$, then there exists an edge $e$ in $H_i$ which contains $v$ and disjoint from the edges of the $(s-1)$-rainbow matching, which produces an $s$-rainbow matching. Hence we may assume the maximum degree in each $H_i$ is at most $3(s-1)$.

    For every color $c$ in $H_i$ and $j\neq i$, let $H_{c}^{j}$ be the subgraph of $H_j$ obtained by deleting all the edges colored by $c$ in $E(H_j)$. Since each $H_j$ is properly colored, there are at most $\frac{n}{2}$ edges colored by $c$ in $E(H_{j})$. So $e(H_{c}^{j})\geq e(H_{j})-\frac{n}{2}>{{n}\choose {2}}-{{n-s+1}\choose {2}}-\frac{n}{2}>{{n}\choose {2}}-{{n-(s-1)+1}\choose {2}}$. By induction, there  exists an $(s-1)$-rainbow matching $\{e_j\}_{j\neq i}$ in $\{H_{v}^{j}\}_{j\neq i}$, which spans $2(s-1)$ vertices $u_{1}, \ldots, u_{2(s-1)}$. Also since  $H_i$ is properly colored, it has at most one edge containing each $u_j$ and colored by $c$. So if the number of edges in $H_i$ colored by $c$ is greater than $2(s-1)$, then there exists an edge $e$ in $H_i$ colored by $c$ and disjoint from $\{e_j\}_{j\neq i}$, which produces an $s$-rainbow matching. So we can now assume that the number of edges in every color in each $H_i$ is at most $2(s-1)$.

    Now pick an arbitrary edge $uv$ in $H_1$. Assume the color of $uv$ is $c(uv)$. Then we delete the vertices $u$, $v$ and all the edges colored by $c(uv)$ in $H_{2}, \ldots, H_{s}$. Denote the resulting graphs by $H'_{2}, \ldots, H'_{s}$. We can see that when $n\geq7s$, for each $i\in \{2, \ldots, s\}$, we have $e(H'_{i})>{{n}\choose {2}}-{{n-s+1}\choose {2}}-2\cdot3(s-1)-2(s-1)>{{n-2}\choose {2}}-{{(n-2)-(s-1)+1}\choose {2}}$. By induction on $s$, there exists an $(s-1)$-rainbow matching in the graphs $H'_{2}, \ldots, H'_{s}$. Taking these $s-1$ edges with the edge $uv$, we obtain an $s$-rainbow matching in $H_{1}$, \ldots, $H_{s}$.

\textbf{Case 2: $k\geq3$.} We do induction on $s$. The case $s=1$ is trivial. We first consider the situation when some $H_i$ has a vertex of degree greater than $k(s-1){{n-2}\choose {k-2}}+s-1$. For every vertex $v\in H_i$ and $j\neq i$, let $H_{v}^{j}$ be the subgraph of $H_j$ induced by the vertex set $V(H_{j})\setminus \{v\}$. Since there are at most ${n-1}\choose {k-1}$ edges containing $v$ in $E(H_{j})$, we have $e(H_{v}^{j})\geq e(H_{j})-{{n-1}\choose{k-1}}>{{n}\choose {k}}-{{n-s+1}\choose {k}}-{{n-1}\choose{k-1}}={{n-1}\choose {k}}-{{(n-1)-(s-1)+1}\choose {k}}$. By induction, there exists an $(s-1)$-rainbow matching $\{e_j\}_{j\neq i}$ in $\{H_{v}^{j}\}_{j\neq i}$, which spans $k(s-1)$ vertices. So if some $H_i$ has a vertex $v$ with degree greater than $k(s-1){{n-2}\choose {k-2}}+s-1$, then there exists an edge $e$ in $E(H_i)$ which contains $v$ and disjoint from the edges of the $(s-1)$-rainbow matching, which produces an $s$-rainbow matching. Hence we may assume the maximum degree in each hypergraph $H_i$ is at most $k(s-1){{n-2}\choose {k-2}}+s-1$.

By induction on $s$, we know that for every $i$ there exists an $(s-1)$-rainbow matching in the hypergraphs $\{H_j\}_{j\neq i}$, spanning $k(s-1)$ vertices. If for some $i$, the $s$-th largest degree of $H_i$ is at most $2(s-1){{n-2}\choose {k-2}}+s-1$, then the sum of degrees of these $k(s-1)$ vertices in $H_i$ is at most \medskip

 $(s-1)[k(s-1){{n-2}\choose {k-2}}+s-1]+(s-1)(k-1)[2(s-1){{n-2}\choose {k-2}}+s-1]=(3k-2)(s-1)^{2}{{n-2}\choose {k-2}}+(s-1)^{2}k$.\medskip

 On the other hand, the maximum degree of the subgraph of $H_i$ by deleting these $k(s-1)$ vertices is at most $s-1$, otherwise, we can find an $s$-rainbow matching. Since $n\geq 3k^{2}s$, we have $e(H_i)>{{n}\choose {k}}-{{n-s+1}\choose {k}}>(s-1)^{2}(3k-\frac{1}{2}){{n-2}\choose {k-2}}>(3k-2)(s-1)^{2}{{n-2}\choose {k-2}}+(s-1)^{2}k+\frac{(s-1)[n-k(s-1)]}{k}$, which guarantees the existence of an edge in $H_i$ disjoint from the previous $(s-1)$-rainbow matching in $\{H_j\}_{j\neq i}$, which produces an $s$-rainbow matching. So we may assume that each $H_i$ contains at least $s$ vertices with degree above  $2(s-1){{n-2}\choose {k-2}}+s-1$.

Now we may greedily select distinct vertices $v_{i}\in V(H_{i})$, such that for each $1\leq i\leq s$, the degree of $v_i$ in $H_i$ exceeds $2(s-1){{n-2}\choose {k-2}}+s-1$. Consider all the subsets of $V(H_{i})\setminus \{v_{1}, \ldots, v_{s}\}$ which together with $v_i$ form an edge of $H_i$. Denote the $(k-1)$-uniform hypergraph by $H_{i}'$. Since each $H_{i}$ is properly colored, we can see that each $H_{i}'$ is rainbow and $e(H_{i}')>2(s-1){{n-2}\choose {k-2}}+s-1-(s-1){{n-2}\choose {k-2}}>{{n-s}\choose {k-1}}-{{n-2s+1}\choose {k-1}}$. By Lemma \ref{lem2}, there exists an $s$-rainbow matching $\{e_{i}\}_{1\leq i\leq s}$ in $\{H_{i}'\}_{1\leq i\leq s}$. Taking the edges $e_{i}\bigcup \{v_i\}$, we obtain an $s$-rainbow matching in $\{H_{i}\}_{1\leq i\leq s}$.

\hfill\vrule height3pt width6pt depth2pt

\section{Concluding Remarks}
In this short note, we propose a generalization of the Erd\H os hypergraph matching conjecture to finding rainbow matchings in  properly-colored hypergraphs, and prove Theorem \ref{thm1} for $s<n/(3k^{2})$. The following conjecture seems plausible.
\begin{conj}
There exists constant $C>0$ such that if $H_{1}$, \ldots, $H_{s}$ are properly-colored $k$-uniform hypergraphs on $n$ vertices, with $n\geq Cks$ and every $e(H_{i})>{{n}\choose {k}}-{{n-s+1}\choose {k}}$, then there exists an $s$-rainbow matching in $H_{1}$, \ldots, $H_{s}$.
\end{conj}
Recall that for the special case when each $H_i$ is identical and rainbow, Frankl and Kupavskii \cite{FK} were able to verify it for $C= 5/3$. However the proof relies on the technique of shifting, while the property of a hypergraph being properly colored may not be preserved under shifting.

It is tempting to believe that Erd\H os Matching Conjecture can be extended to properly-colored hypergraphs for the entire range of $s$, that is, once the number of edges in each hypergraph exceeds the maximum of $\binom{n}{k}-\binom{n-s+1}{k}$ and $\binom{ks-1}{k}$, then one can find an $s$-rainbow matching. However this is false in general, a simple construction is by taking $s=2$ and $n=2k$. The maximum of these two expressions is $\binom{2k-1}{k}$, while one can let $H_1$ be a rainbow $K_{2k}^k$ with an edge coloring $c_1$, and $H_2$ be on the same vertex set with edge coloring $c_2$, such that $c_2(e)=c_1([2k]\setminus e)$. Then clearly each $H_i$ contains $\binom{2k}{k}>\binom{2k-1}{k}$ edges and there is no $2$-rainbow matching. It would be interesting to find constructions for $s$ close to $n/k$, and formulate a complete conjecture for properly-colored hypergraphs.

%%%%%%%
%%%%%%%
%%%%%%%
%%%%%%%
%%%%%%%

%\newpage

\end{document}